\documentclass[3p,times]{elsarticle}

\usepackage{enumitem}
\usepackage{hyperref}
\usepackage{amsmath}
\usepackage{pifont}
\usepackage[utf8]{inputenc}
\makeatletter
\def\LaTeX{\leavevmode L\raise.42ex
    \hbox{\kern-.3em\size{\sf@size}{0pt}\selectfont A}\kern-.15em\TeX}
\makeatother

\newcommand{\BibTeX}{{\rm B\kern-.05em{\sc
          i\kern-.025emb}\kern-.08em\TeX}}

\makeatletter
\def\@currentlabel{2.1}\label{e:dispaa}
\def\@currentlabel{2.21}\label{e:dispau}
\def\@currentlabel{2.22}\label{e:dispav}
\def\@currentlabel{2.23}\label{e:dispaw}
\def\@currentlabel{2.24}\label{e:dispax}
\def\theequation{\thesection.\@arabic\c@equation}
\makeatother

\renewcommand{\theequation}{\arabic{section}.\arabic{equation}}

\newcommand{\R}{\mathbb R}

\def \O{\Omega}
\everymath{\displaystyle}
\newtheorem{theorem}{Theorem}
\renewcommand*{\thetheorem}{\Alph{theorem}}
\newtheorem{theo}{Theorem} [section]
\newtheorem{lem}{Lemma} [section]
\newtheorem{prop}{Proposition} [section]

\newtheorem{definition}{Definition} [section]

\newtheorem{rem}{Remark}[section]

\renewcommand{\theequation}{\thesection.\arabic{equation}}
\renewcommand{\thesection}{\arabic{section}}
\renewcommand{\theequation}{\thesection.\arabic{equation}}
\let\ssection=\section\renewcommand{\section}{\setcounter{equation}{0}\ssection}
\begin{document}
\begin{frontmatter}
\title{Existence and multiplicity results for a new $p(x)$-Kirchhoff problem}
\author[mk1,mk2]{Mohamed Karim Hamdani}
\ead{hamdanikarim42@gmail.com}
\author[ah3,ah1,ah2]{ Abdellaziz Harrabi}
\ead{abdellaziz.harrabi@yahoo.fr}
\author[Mtiri]{ Foued Mtiri}
\ead{mtirifoued@yahoo.fr}
\author[dr1,dr2,dr3]{ Du\v{s}an D. Repov\v{s}\corref{cor1}}
\ead{dusan.repovs@guest.arnes.si}
\begin{center}
\address[mk1]{ Mathematics Department, University of Sfax, Faculty of Science of Sfax, Sfax, Tunisia.}
\address[mk2] {Military School of Aeronautical Specialities, Sfax, Tunisia.}
\address[ah3]{Mathematics Department, University of Kairouan, Higher Institute of Applied Mathematics and Informatics, Kairouan, Tunisia.}
\address[ah1]{ Mathematics Department, Northern Borders University, Arar, Saudi Arabia.}
\address[ah2]{Senior associate in the Abdus Salam International Centre for Theoretical Physics, Trieste, Italy.}
\address[Mtiri]{ Mathematics Department, University of King Khaled, Faculty of Sciences and Arts of Muhail Asir, Abha Saudi Arabia.}
\address[dr1] {Faculty of Education, University of Ljubljana, Ljubljana, Slovenia}
\address[dr2] {Faculty of Mathematics and Physics, University of Ljubljana, Ljubljana, Slovenia}
\address[dr3] {Institute of Mathematics, Physics and Mechanics,  Ljubljana, Slovenia}
\cortext[cor1]{Corresponding author: Du\v{s}an D. Repov\v{s}, dusan.repovs@guest.arnes.si}

\end{center}\begin{abstract}
In this work, we study the existence and multiplicity results for the following nonlocal $p(x)$-Kirchhoff problem:
\begin{equation}
 \label{10}
\begin{cases}
-\left(a-b\int_\Omega\frac{1}{p(x)}| \nabla u| ^{p(x)}dx\right)div(|\nabla u| ^{p(x)-2}\nabla u)=\lambda |u| ^{p(x)-2}u+g(x,u) \mbox{ in } \Omega, \\
u=0,\mbox{ on } \partial\Omega,
\end{cases}
\end{equation}
where $a\geq b > 0$ are constants, $\Omega\subset \mathbb{R}^N$ is a bounded smooth domain, $p\in C(\overline{\Omega})$ with $N>p(x)>1$, $\lambda$ is a real parameter and $g$ is a continuous function. The analysis developed in this paper proposes an approach based on the idea of considering a new nonlocal term which presents interesting difficulties.
\end{abstract}

\end{frontmatter}
\section{Introduction and statement of the main results}
\label{sect1}
In this work, we study the existence and multiplicity results for the following nonlocal $p(x)$-Kirchhoff problem:
\begin{equation}
 \label{prob1}
\begin{cases}
-\left(a-b\int_\Omega\frac{1}{p(x)}| \nabla u| ^{p(x)}dx\right)div(| \nabla u| ^{p(x)-2}\nabla u)=\lambda | u| ^{p(x)-2}u+g(x,u) \mbox{ in } \Omega, \\
u=0,\mbox{ on } \partial\Omega,
\end{cases}
\end{equation}
where $\Omega\subset \R^N$ is a bounded smooth domain, $p\in C(\overline{\O})$ with $N>p(x)>1$, $a, b > 0$ are constants, $g$ is a continuous function satisfying conditions which will be stated later, $\lambda>0$ is a real parameter and $div(| \nabla u| ^{p(x)-2}\nabla u)$ is the $p(x)$-Laplacian operator, that is,
\[
\Delta_{p(x)}=div(| \nabla u| ^{p(x)-2}\nabla u)=\sum_{i=1}^{N}\left(| \nabla u| ^{p(x)-2}\frac{\partial u}{\partial x_i}\right),
\]
which is not homogeneous and is related to the variable exponent Lebesgue space $L^{p(x)}(\O)$ and the variable exponent Sobolev space $W^{1,p(x)}(\O)$.
These facts imply some difficulties. For example, some classical theories and methods, including the Lagrange multiplier theorem and the theory of Sobolev spaces, cannot be
applied. We call \eqref{prob1} a problem of Kirchhoff type because of the appearance of the term
 \[
 b\int_\O\frac{1}{p(x)}| \nabla u| ^{p(x)}dx
 \]
 which makes the study of \eqref{prob1} interesting.

In the previous decades, the Kirchhoff type problem \eqref{prob1} with $p(x)\equiv 2$  has been
the object of intensive research due to its strong relevance in applications (see~\cite{LL,LY2,WSL}). Indeed, the study of
Kirchhoff type problems, which arise in various models of physical and biological systems, has received more and more attention in recent years.
More precisely, Kirchhoff  established a model given
by the equation\begin{equation} \label{kirchhoff}
\rho\frac{\partial ^2u}{\partial
t^2}-\Big(\frac{p_0}{h}+\frac{E}{2L} \int_0^L\left|\frac{\partial
u}{\partial x}\right|^2dx\Big) \frac{\partial ^2u}{\partial
x^2}=0,
\end{equation}
where $\rho$, $p_0$, $h$, $E$, $L$ are constants which represent
some physical meanings respectively. Eq.~\eqref{kirchhoff} extends
the classical D'Alembert wave equation by considering the effects of
the changes in the length of the strings during the vibrations.

Since the variable exponent spaces  have been thoroughly
studied by  Kov\'a\u{c}ik and  R\'akosn\'ik~\cite{KR},
they have been used in the previous decades to model various phenomena.
In the studies of a class of non-standard variational problems and PDEs,
variable exponent spaces play an important role for example, in electrorheological
fluids~\cite{R,RR1996,RR2001}, thermorheological fluids~\cite{AR},
image processing~\cite{AMS,CLR,LLP}, etc.
In recent years, there has been
a great deal of work done on problem \eqref{prob1}, especially concerning the
existence, multiplicity, uniqueness and regularity of solutions. Some important
and interesting results can be found, for example,
in~\cite{AB,ACX,AH2,AM2,BO,C,CF,CP,DH,DM,F,FZ1,FZ2,FZ3,H,HH,KR2,KR,MR,MR2007,PRD,PZC,RD,ZZZ}
and
references therein.

At first, the eigenvalues of $p(x)$-Laplacian Dirichlet problem were
studied in~\cite{FZZ}, i.e.,~if $\O\subset \R^N$ is a smooth bounded
domain, the Rayleigh quotient
\begin{equation}
\label{eigen}
\lambda_{p(.)}=\inf_{u\in W_0^{1,p(x)}(\O)\setminus \{0\}}\frac{\int_\O\frac{1}{p(x)}| \nabla u| ^{p(x)}dx}{\int_\O\frac{1}{p(x)}| u| ^{p(x)}dx}
\end{equation}
is zero in general, and only under some special conditions $\lambda_{p(.)}> 0$ holds. For example, when $\O\subset \R \ (N=1)$ is an interval,
results show that $\lambda_{p(.)}> 0$ if and only if $p(.)$ is monotone. It is well known that $\lambda_p > 0$ plays a very important role in the study
of $p$-Laplacian problems.

Motivated by the papers mentioned above, our main purpose is to consider the perturbed problem \eqref{prob1} with a new nonlocal term
\[
a-b\int_\O\frac{1}{p(x)}| \nabla u| ^{p(x)}dx
\]
which presents interesting difficulties.
The key argument in our main result is the proof that the energy functional $J$ (which appeared in \eqref{funct energ}) of problem \eqref{prob1} possesses a Mountain Pass energy $c$.

To deal with the difficulty caused by the noncompactness due to the Kirchhoff function term, we must estimate precisely the value of $c$ and give a threshold value (see Lemma \ref{lemme1}) under which the Palais--Smale condition at the level $c$ for $J$ is satisfied. So the variational technique for problem \eqref{prob1} becomes more delicate. We  obtain a nontrivial weak solution by using the Mountain Pass theorem. To the best of our knowledge, the present papers results have not been covered yet in the literature.

Suppose that the nonlinearity $g (x, t)\in  C(\overline{\O}\times \R)$ satisfies the following assumptions:
\begin{enumerate}[label={$g\sb{\arabic*}$:},ref={{$g\sb{\arabic*}$}}]
\item\label{hypot:g1}  the subcritical growth condition:
\[
|g(x,s)| \leq C(1+|s|^{q(x)-1}),  \mbox{ for all } (x,s)\in \O\times \R,
\]
where $C>0$ and $p(x)<q(x)<p^*(x); $
\item\label{hypot:g2} $ \lim_{s\rightarrow 0} \frac{g(x,s)}{|s|^{p(x)-2}s}=0; $
\item\label{hypot:g3} there exist $s_A > 0$ and $\theta\in (p^+, \frac{2p^{-^2}}{p^+})$ such that
\[
 0<\theta G(x,s)\leq sg(x,s),\;\mbox{for all} \; |s| \geq s_A,\; x\in \O,
\]
where $G(x,s)=\int_{0}^{s}g(x,t)dt$;
\item\label{hypot:g4} $g(x,-s)=-g(x,s) \mbox{ for all } (x,s)\in \O\times \R$.
\end{enumerate}
 Now we can state our main results:
\begin{theo}
\label{theo1.1}
Suppose that the function $q\in C(\overline{\O})$ satisfies
\begin{equation}
\label{cond q}
1<p^-<p(x)<p^+<2p^-<q^-<q(x)<p^*(x).
\end{equation}
Then for any $\lambda\in\R$, with \eqref{hypot:g1}--\eqref{hypot:g3} satisfied,
problem
\eqref{prob1} has a nontrivial weak solution.
\end{theo}
\begin{theo}
\label{multiplicity results}
Suppose that the function $q\in C(\overline{\O})$ satisfies
\[
\label{cond q2}
1<p^-<p(x)<p^+<2p^-<q^-<q(x)<p^*(x).
\]
 Then for any $\lambda\in\R$, with \eqref{hypot:g1}--\eqref{hypot:g4} satisfied,
 problem
 \eqref{prob1} has infinitely many solutions $\{u_n\}$ such that $I(u_n)\to\infty$  as $n\to\infty$.
\end{theo}
\begin{rem}
Hypothesis \eqref{hypot:g3} is known as the Ambrosetti--Rabinowitz's
superlinear condition (see~\cite{CP}). Moreover, condition
\eqref{hypot:g3} ensures that the Euler--Lagrange functional associated with
 problem \eqref{prob1} possesses the geometry of
Mountain Pass theorem and it also guarantees the boundedness of the Palais--Smale sequence corresponding to the Euler--Lagrange's functional.
\end{rem}
This paper is organized as follows.
In Section \ref{section2}, we present some necessary preliminary knowledge on variable exponent Sobolev spaces.
In Section \ref{section3}, we prove  the Palais-Smale compactness condition.
In Section \ref{section4}, we prove  Theorem \ref{theo1.1} via the Mountain Pass theorem.
In Section \ref{section5}, we prove  Theorem  \ref{multiplicity results} via the Fountain theorem.
In this paper, $|\cdot| $ denotes the Lebesgue measure on $\O,$ and
$C$ (respectively, $C_{\epsilon}$)  always denotes a generic positive constant
independent of $n$ and $\epsilon$ (respectively, $n$).

\section{Preliminaries on variable exponent  spaces}\label{section2}
In order to discuss problem \eqref{prob1}, we need some theory on spaces $L^{p(x)}(\O)$ and $W^{1,p(x)}(\O)$ which
we shall call generalized Lebesgue Sobolev spaces. Let $\O$
 be a bounded domain of $\R^N$, denote $C_+(\overline{\O})=\{p(x); p(x)\in C(\overline{\O}), p(x)>1, \;\mbox{for all}  \; x\in \overline{\O}\}$ and
$p^-=\inf_\O p(x)\leq p(x)\leq p^+=\sup_\O p(x)<N$.\\
For any $p\in C_+(\overline{\O})$, we introduce the variable exponent Lebesgue space
$$L^{p(\cdot)}(\O)=\left\{u: u \mbox{ is a measurable real-valued function such that } \int_\O |u(x)|^{p(x)}dx<\infty\right\},$$
endowed with the so-called Luxemburg norm

$$\|u\|_{L^{p(x)}(\O)}=|u|_{p(.)}=\inf \left\{\mu>0;\int_\O\left|\frac{u(x)}{\mu}\right|^{p(x)}dx\leq 1\right\},$$
which is a separable and reflexive Banach space. For basic properties of the variable exponent Lebesgue spaces we
refer to \cite{FF,KR,Y}.

\begin{prop}[\cite{Y}]
 The space $(L^{
p(x)}
(\O), |.|_{p(x)})$ is separable, uniformly convex,
reflexive, and its conjugate space is $(L^{
q(x)}
(\O), |.|_{q(x)}),$ where $q(x)$ is the conjugate function
of $p(x)$ i.e
\[
\frac{1}{p(x)}+\frac{1}{q(x)}=1,\;\;\mbox{for all} \; x\in \O.
\]
For all $u\in L^{p(x)}(\O), v\in L^{q(x)}(\O),$ the H\"older type inequality
\[
  \Big|\int_\O uvdx\Big|  \leq\left(\frac{1}{p^-}+\frac{1}{q^-}\right)|u| _{p(x)}|v|_{q(x)}
\]
  holds.
\end{prop}
The inclusion between Lebesgue spaces also generalizes the classical framework,
namely, if $0<|\O| <\infty$ and $p_1$, $p_2$ are variable exponents such that $p_1 \leq p_2$ in $\O,$ then there exists a continuous embedding $L^{p_2(x)}(\O)\to L^{p_1(x)}(\O)$.

An important role in working with the generalized Lebesgue--Sobolev spaces is played by the $m(\cdot)$-modular of the
$L^{p(\cdot)}(\O)$ space, which is the modular $\rho_{p(\cdot)}$ of the space $L^{p(\cdot)}(\O)$
\[
   \rho_{p(\cdot)}(u):=\int_{\Omega} |u|^{p(x)} \,dx.
\]
\begin{lem}[\cite{DHHR}]
\label{lemmaineq}
Suppose that $u_n, u \in L^{p(\cdot)}$ and $p_{+} < +\infty$. Then the following
properties hold:
\begin{enumerate}
\item $| u| _{p(\cdot)} > 1 \Rightarrow
| u| _{p(\cdot)}^{p_{-}} \leq \rho_{p(\cdot)}(u) \leq | u| _{p(\cdot)}^{p^{+}}$;

\item $| u| _{p(\cdot)} < 1 \Rightarrow   | u| _{p(\cdot)}^{p^{+}}
\leq \rho_{p(\cdot)}(u) \leq | u| _{p(\cdot)}^{p^{-}}$;

\item $| u| _{p(\cdot)} < 1$ (respectively, $= 1; {>} 1) \Longleftrightarrow
 \rho_{p(\cdot)}(u) < 1$ (respectively, ${=} 1; {>} 1$);

\item $| u_n| _{p(\cdot)} \to 0$ (respectively,
$\to +\infty) \Longleftrightarrow \rho_{p(\cdot)}(u_n) \to 0$
 (respectively, $\to +\infty$);

\item $\lim_{n\to \infty}| u_n-u| _{p(x)}=0 \Longleftrightarrow \lim_{n\to \infty}\rho_{p(\cdot)}(u_n-u)=0$.
\end{enumerate}
\end{lem}
The Sobolev space with variable exponent $W^{1,p(x)}(\O)$ is defined as
\begin{equation*}
   W^{1,p(x)}(\Omega):=\Big\{u: \O\subset \R^N\to \R: u \in L^{p(x)}(\Omega), |\nabla u| \in L^{p(x)}(\Omega)\Big\},
\end{equation*}
equipped with the norm
$$\|u\|_{1,p(x)}=\|u\|_{p(x)}+\|\nabla u\|_{p(x)}.$$

Then $W_0^{1,p(x)}(\O)$ is defined as the closure of $C_0^\infty (\O)$ with respect to
the norm $\| u\| _{1,p(x)}$. In this way, $L^{p(x)}(\Omega)$, $W_0^{1,p(x)}(\O)$ and $W^{1,p(x)}(\O)$
become separable and reflexive Banach spaces. For more details, we refer
to~\cite{DHHR,FSZ,FZ1}. Moreover, we define
\[
  p^*(x) =
\begin{cases}
   \frac{Np(x)}{N-p(x)}, & \mbox{ if } p(x)<N \\
    +\infty, & \mbox{ if } p(x)\geq N.
\end{cases}
\]
The following results were proved in~\cite{FZ1}.
\begin{prop}[{Sobolev Embedding \cite{FZ1}}]
\label{sobolev} For $p,q\in C_+(\overline{\O})$ such that $1\leq q(x)\leq p^*(x)$ for all $x\in \overline{\O}$, there is a continuous embedding
\[
W^{1,p(x)}(\O)\hookrightarrow L^{q(x)}(\O).
\]
If we replace $\leq$ with $<$, the embedding is compact.
\end{prop}
\begin{prop}[{Poincar\'e Inequality \cite{FZ1}}]
\label{poinc} There is a constant $C > 0$, such that
\begin{equation}
\label{point}
\| u\| _{L^{p(x)}(\O)}\leq C \| \nabla u\| _{L^{p(x)}(\O)}
\end{equation}
for all $u\in W_0^{1,p(x)}(\O)$.
\end{prop}
\begin{rem}
By Proposition \ref{poinc}, we know that $\| \nabla u\| _{L^{p(x)}(\O)}$ and $\| u\| _{W^{1,p(x)}(\O)}$ are equivalent
norms on $W_0^{1,p(x)}(\O)$.
\end{rem}

\begin{lem}[\cite{FZ2}]\label{s+}
Denote
\[
A(u)=\int_\O\frac{1}{p(x)}|\nabla u| ^{p(x)}dx,\;\mbox{ for all }u\in W_0^{1,p(x)}(\O).
\]
Then $A(u)\in C^1(W_0^{1,p(x)}(\O),\R)$  and the derivative operator $A'$ of $A$ is
\[
\langle A'(u),v\rangle=\int_\O|\nabla u|^{p(x)-2}\nabla u \nabla vdx \mbox{ for all }u,v\in W_0^{1,p(x)}(\O),
\]
and the following holds:
\begin{enumerate}
  \item $A$ is a convex functional;
  \item $A': W_0^{1,p(x)}(\O)\to (W^{-1,p'(x)}(\O))=\left(W_0^{1,p(x)}(\O)\right)^*$ is a bounded homeomorphism and strictly monotone operator, and the conjugate exponent satisfies $\frac{1}{p(x)}+\frac{1}{p'(x)}=1$;
  \item $A'$ is a mapping of type $S_+$, namely, $u_n\rightharpoonup u$ and $\limsup \langle A'(u_n),u_n-u\rangle\leq 0$,
imply $u_n\to u$ (strongly) in $W_0^{1,p(x)}(\O)$.
\end{enumerate}
\end{lem}
\begin{definition}
 We say that $u\in W_0^{1,p(x)}(\O)$ is a weak solution of \eqref{prob1}, if
\[
\left(a-b\int_\O\frac{1}{p(x)}|\nabla u| ^{p(x)}dx\right)\int_\O|\nabla u|^{p(x)-2}\nabla u \nabla \varphi dx-\lambda\int_{\O}|u|^{p(x)-2}u \varphi dx= \int_{\O} g(x,u)\varphi dx,
\]
where $\varphi\in W_0^{1,p(x)}(\O)$.
\end{definition}
The energy functional $J: W_0^{1,p(x)}(\O) \to \mathbb{R}$ associated with problem \eqref{prob1}
\begin{equation}
\label{funct energ}
J(u)=a\int_\O\frac{1}{p(x)}|\nabla u| ^{p(x)}dx-\frac b2\left(\int_\O\frac{1}{p(x)}| \nabla u| ^{p(x)}dx\right)^2-\lambda\int_{\O}\frac{1}{p(x)}|  u| ^{p(x)}dx- \int_{\O} G(x,u)dx,
\end{equation}
for all $u\in W_0^{1,p(x)}(\O)$ is well defined and of $C^1$ class on $W_0^{1,p(x)}(\O)$. Moreover, we have
\begin{eqnarray}\label{deriv funct energ}
\langle J'(u),\varphi\rangle&=&\left(a-b\int_\O\frac{1}{p(x)}|\nabla u|^{p(x)}dx\right)\int_\O|\nabla u|^{p(x)-2}\nabla u \nabla \varphi dx-\lambda\int_{\O}|u|^{p(x)-2}u \varphi dx\nonumber\\&\quad-& \int_{\O} g(x,u)\varphi dx,
\end{eqnarray}
for all $u, \varphi \in W_0^{1,p(x)}(\O)$.
Hence, we can observe that the critical points of the functional $J$ are the weak
solutions for problem \eqref{prob1}. In order to simplify the presentation we will denote the norm of $W_0^{1,p(x)}(\O)$ by $\| .\|$, instead of $\| \cdot\| _{W_0^{1,p(x)}(\O)}$.
\section{The Palais-Smale Compactness Condition}\label{section3}
Recall now the definition of the Palais-Smale compactness condition.
\begin{definition}
\label{def1} Let $(W_0^{1,p(x)}(\O),\;\|.\| )$ be a Banach space and $J \in C^1(W_0^{1,p(x)}(\O))$. Given $c \in\R$, we say that $J$ satisfies the Palais--Smale condition at the level $c\in \R$ (``$(PS)_c$ condition'', for short) if every sequence $\{u_n\} \in W_0^{1,p(x)}(\O)$ satisfying
\begin{equation}
\label{condps}
 J(u_n) \rightarrow c \mbox{ and } J'(u_n)\rightarrow  0  \mbox{ in } W^{-1,p'(x)}(\O) \mbox{ as } \;n\rightarrow \infty,
\end{equation}
has a convergent subsequence.
\end{definition}
First, we investigate the compactness conditions for the functional $J$.
\begin{lem}
 \label{lemme1}
Assume that \eqref{hypot:g1}--\eqref{hypot:g3} hold. Then the functional $J$ satisfies the $(PS)_c$ condition, where precisely $c<\frac{a^2}{2b}$.
\end{lem}
{\bf Proof.} We proceed in two steps.\\
\textbf{Step 1.} We prove that $\{u_n\}$ is bounded in $W_0^{1,p(x)}(\O)$. Let $\{u_n\} \subset W_0^{1,p(x)}(\O)$ be a $(PS)_c$ sequence such that $c<\frac{a^2}{2b}$.
\\$\mathbf{\bullet}$ For $\lambda\leq0$. From \eqref{condps} and \eqref{hypot:g3}, for n large enough, we have
\begin{align*}
\label{ing2}
C&+\| u_n\| \\&\geq\theta J(u_n)-\langle J'(u_n), {u_n}\rangle\\&\geq \theta \left(a\int_\O\frac{1}{p(x)}| \nabla u_n| ^{p(x)}dx-\frac b2\left(\int_\O\frac{1}{p(x)}| \nabla u_n| ^{p(x)}dx\right)^2-\lambda\int_{\O}\frac{1}{p(x)}|  u_n| ^{p(x)}dx- \int_{\O} G(x,u_n)dx \right)\nonumber\\
&- \left(\left[a-b\int_\O\frac{1}{p(x)}| \nabla u_n| ^{p(x)}dx\right]\int_\O| \nabla u_n| ^{p(x)}dx-\lambda\int_{\O}| u_n| ^{p(x)}dx- \int_{\O} g(x,u_n)u_ndx \right)\nonumber\\&\geq a(\frac{\theta}{p^+}-1)\int_\O| \nabla u_n| ^{p(x)}dx+b(\frac{-\theta}{2p^{-^2}}+\frac{1}{p^+})\left(\int_\O| \nabla u_n| ^{p(x)}dx\right)^2-\lambda(\frac{\theta}{p^-}-1)\int_{\O}| u_n| ^{p(x)}dx-C | \O| ,\nonumber
\end{align*}
where $| \O| =\int_\O dx$. Since $\lambda\leq0$,  we can deduce that
\[
\label{ing2}
C+\| u_n\| \geq a(\frac{\theta}{p^+}-1)\| u_n\| ^{p^-}+b(\frac{-\theta}{2p^{-^2}}+\frac{1}{p^+})\| u_n\| ^{2p^-}-C | \O| .
\]
It follows from  \eqref{cond q} that $\{u_n\}$ is
bounded in $W_0^{1,p(x)}(\O)$.
\\$\mathbf{\bullet}$ For $\lambda>0$. Arguing by contradiction, we assume that,
passing eventually to a subsequence, still denoted by $\{u_n\}$, we have $\| u_n\| \to +\infty$ as $n\to +\infty$. By \eqref{condps} and \eqref{hypot:g3}, for n large enough, we have
\begin{eqnarray*}
\label{ing2}
C&+&\| u_n\| \geq\theta J(u_n)-\langle J'(u_n), {u_n}\rangle\\&\geq& a(\frac{\theta}{p^+}-1)\int_\O| \nabla u_n| ^{p(x)}dx+b(\frac{-\theta}{2p^{-^2}}+\frac{1}{p^+})\left(\int_\O| \nabla u_n| ^{p(x)}dx\right)^2-\lambda(\frac{\theta}{p^-}-1)\int_{\O}| u_n| ^{p(x)}dx-C | \O| ,\nonumber
\end{eqnarray*}
 Therefore the last inequality, together with \eqref{point}, implies that
\[
\label{ing2}
C+\| u_n\| +\lambda C(\frac{\theta}{p^-}-1)\| u_n\| ^{p^+}\geq a(\frac{\theta}{p^+}-1)\| u_n\| ^{p^-}+b(\frac{-\theta}{2p^{-^2}}+\frac{1}{p^+})\| u_n\| ^{2p^-}-C | \O| .
\]
Dividing the above inequality by $\| u_n\| ^{p^+}$, taking into account \eqref{cond q} holds and
passing to the limit as $n\to \infty$, we obtain a contradiction. It follows that $\{u_n\}$ is
bounded in $W_0^{1,p(x)}(\O)$.
\\\textbf{Step 2.}
Here, we will prove that $\{u_n\}$ has
a convergent subsequence in $W_0^{1,p(x)}(\O)$. It follows from  Proposition \ref{sobolev}
that the embedding
\[
W_0^{1,p(x)}(\O) \hookrightarrow L^{s(x)}(\O)
\]
 is compact, where $1\leq s(x)<p(x)^*$. Passing, if necessary, to a subsequence, there exists $u\in W_0^{1,p(x)}(\O)$ such that
\begin{equation}
\label{cvg}
u_n\rightharpoonup u\mbox{ in } W_0^{1,p(x)}(\O),\; u_n \to u \mbox{ in } L^{s(x)}(\O),\;\; u_n(x)\to u(x), \mbox{ a.e. in } \O.
\end{equation}
By H\"older's inequality and \eqref{cvg}, we have
\begin{eqnarray*}
{  \left|\int_{\O}{| u_n|} ^{p(x)-2}u_n (u_n-u)dx\right|}   &\leq& \int_{\O}{| u_n| }^{p(x)-1}| u_n-u| dx\\&\leq& {\Big|{|  u_n|} ^{p(x)-1}\Big|} _{\frac{p(x)}{p(x)-1}}| u_n-u| _{p(x)}\\&\to& 0 \mbox{ as } n\to \infty
\end{eqnarray*}
and thus,
\begin{equation}
\label{cvgu}
\lim_{n\to \infty}\int_{\O}| u_n| ^{p(x)-2}u_n (u_n-u)dx=0.
\end{equation}
 By virtue of conditions \eqref{hypot:g1} and  \eqref{hypot:g2}, one has that for every $\epsilon \in (0,1)$ there exists $C_{\epsilon}>0$ such that
\begin{equation}
\label{ss}
 | g(x,u_n)| \leq \epsilon | u_n| ^{p(x)-1}+C_\epsilon | u_n| ^{q(x)-1}.
\end{equation}
By \eqref{ss} and Proposition \ref{sobolev}, it follows that
\begin{eqnarray*}
\Big|\int_\O g(x,u_n)(u_n-u)dx\Big|  &\leq& \int_\O \epsilon | u_n| ^{p(x)-1}| u_n-u| +C_\epsilon | u_n| ^{q(x)-1}| u_n-u| dx\\&\leq& \epsilon{\Big|{|  u_n|} ^{p(x)-1}\Big|}  _{\frac{p(x)}{p(x)-1}}{| u_n-u|} _{p(x)}+C_\epsilon \epsilon{\Big|{|u_n|} ^{q(x)-1}\Big|} _{\frac{q(x)}{q(x)-1}}| u_n-u| _{q(x)}\\&\to& 0 \mbox{ as } n\to \infty
\end{eqnarray*}
which shows that
\begin{equation}
\label{cvgf}
\lim_{n\to \infty}\int_{\O}g(x,u_n)(u_n-u)dx=0.
\end{equation}
By \eqref{condps}, we have
\[
\langle J'(u_n),u_n-u\rangle \to 0.
\]
 Therefore
\begin{eqnarray*}
\label{deriv funct energ}
\langle J'(u_n),u_n-u\rangle&=\left(a-b\int_\O\frac{1}{p(x)}| \nabla u_n| ^{p(x)}dx\right)\int_\O| \nabla u_n| ^{p(x)-2}\nabla u_n (\nabla u_n-\nabla u)dx\\
&-\lambda\int_{\O}| u_n| ^{p(x)-2}u_n (u_n-u)dx- \int_{\O} g(x,u_n)(u_n-u) dx\to 0.
\end{eqnarray*}
So, we can deduce from \eqref{cvgu} and \eqref{cvgf} that
\begin{equation}
\label{gj}
\left(a-b\int_\O\frac{1}{p(x)}|\nabla u_n|^{p(x)}dx\right)\int_\O|\nabla u_n|^{p(x)-2}\nabla u_n (\nabla u_n-\nabla u)dx\to 0.
\end{equation}
Since $\{u_n\}$ is bounded in $W_0^{1,p(x)}(\O)$, passing to a subsequence, if necessary, we may assume that
\[
\int_\O\frac{1}{p(x)}|\nabla u_n|^{p(x)}dx\to t_0\geq 0\mbox{ as } n\to \infty.
\]
{\bf Case 1.} If $t_0 = 0$ then $\{u_n\}$ strongly converges to $u = 0$ in $W_0^{1,p(x)}(\O)$ and the proof is finished.
\\
{\bf Case 2.} If $t_0 > 0$ we need to consider two subcases:
\\
{\bf Subcase 1.} If $t_0\neq \frac ab$ then $a-b\int_\O\frac{1}{p(x)}|\nabla u_n|^{p(x)}dx\to 0$ is not true and no
subsequence of $\{a-b\int_\O\frac{1}{p(x)}|\nabla u_n|^{p(x)}dx\to 0\}$  converges to zero. Therefore, there exists $\delta > 0$
such that $\left| a-b\int_\O\frac{1}{p(x)}|\nabla u_n|^{p(x)}dx\right|  >\delta>0$ when $n$ is large enough. So, it is clear that
\begin{equation}
\{a-b\int_\O\frac{1}{p(x)}|\nabla u_n|^{p(x)}dx\to 0\}  \mbox{ is
bounded}.
\end{equation}
{\bf Subcase 2.}\footnote{This case does not exist if the Kirchhoff function is given by $a+b\int_\O\frac{1}{p(x)}|\nabla u_n|^{p(x)}dx$.} If $t_0 =\frac ab$ then $a-b\int_\O\frac{1}{p(x)}|\nabla u_n|^{p(x)}dx\to 0$.

 We define
\[
\varphi(u)=\lambda\int_\O\frac{1}{p(x)}|u|^{p(x)}dx+\int_\O G(x,u)dx,\; \mbox{for all} \; u\in W_0^{1,p(x)}(\O).
\]
Then
\[
\langle \varphi'(u),v\rangle=\lambda\int_\O|u|^{p(x)-2}u  vdx +\int_\O g(x,u)vdx,\; \mbox{for all} \; u,v\in W_0^{1,p(x)}(\O).
\]
It follows that
\[
\langle \varphi'(u_n)-\varphi'(u),v\rangle=\lambda\int_\O(|u_n|^{p(x)-2}u_n-|u|^{p(x)-2}u)vdx
+\int_\O (g(x,u_n)-g(x,u))vdx.
\]
To complete the argument we need the following lemma.
\begin{lem}
\label{jkjk}
Let $u_n,u\in W_0^{1,p(x)}(\O)$  such that \eqref{cvg} holds. Then, passing to a subsequence, if necessary, the following properties hold:
\begin{enumerate}
  \item[(i)] $\int_\O(|u_n|^{p(x)-2}u_n-|u|^{p(x)-2}u)vdx=0$;
  \item[(ii)] $\lim_{n\to \infty}\int_\O|g(x,u_n)-g(x,u)||v|dx=0$;
  \item[(iii)] $\langle \varphi'(u_n)-\varphi'(u),v\rangle\to 0,\;v\in W_0^{1,p(x)}(\O)$.
\end{enumerate}
\end{lem}
{\bf Proof.} By \eqref{cvg}, we have  $u_n\to u$ in $L^{p(x)}(\O)$ which implies that
\begin{equation}
\label{cvgi}
|u_n|^{p(x)-2}u_n\to |u|^{p(x)-2}u \mbox{ in } L^{\frac{p(x)}{p(x)-1}}(\O).
\end{equation}
Due to H\"older's inequality, we have
\begin{eqnarray}
\left|\int_\O(|u_n|^{p(x)-2}u_n-|u|^{p(x)-2}u)vdx\right|&\leq& \int_\O||u_n|^{p(x)-2}u_n-|u|^{p(x)-2}u||v|dx\nonumber\\
&\leq& \Big||u_n|^{p(x)-2}u_n-|u|^{p(x)-2}u\Big|_{\frac{p(x)}{p(x)-1}}|v|_{p(x)}\nonumber\\
&\leq& C \Big||u_n|^{p(x)-2}u_n-|u|^{p(x)-2}u\Big|_{\frac{p(x)}{p(x)-1}}\|v\|\nonumber\\
&\to& 0.
\end{eqnarray}
 By a slight modification of the proof above, we
can also prove part $(ii)$ so we omit the details.
\[\int_\O|g(x,u_n)-g(x,u)||v|dx\leq \int_\O [\epsilon(|u_n|^{p(x)-2}u_n-|u|^{p(x)-2}u)+C_\epsilon (|u_n|^{q(x)-1}-|u|^{q(x)-1})]|v|dx \to0.
\]
Finally, part $(iii)$ follows by combining parts $(i)$ and $(ii)$.
Consequently, $\|  \varphi'(u_n)-\varphi'(u)\| _{W^{-1,p'(x)}}\to 0$ and $\varphi'(u_n)\to\varphi'(u)$.

We can now complete the proof of Subcase 2:

 By Lemma \ref{jkjk} and since
 $\langle J'(u),u\rangle=\left(a-b\int_\O\frac{1}{p(x)}|\nabla u|^{p(x)}dx\right)\int_\O|\nabla u|^{p(x)-2}\nabla u \nabla \varphi dx-\langle \varphi'(u),v\rangle$, $\langle J'(u),u\rangle\to 0$ and $a-b\int_\O\frac{1}{p(x)}|\nabla u|^{p(x)}dx\to 0$, it follows that $\varphi'(u_n)\to 0\;(n\to\infty)$, i.e.,~\[
\langle \varphi'(u),v\rangle=\lambda\int_\O|u|^{p(x)-2}u  vdx +\int_\O g(x,u)vdx,\; \mbox{for all} \; v\in W_0^{1,p(x)}(\O),
\]
 and therefore
\[
\lambda |u(x)|^{p(x)-2}u(x)+g(x,u(x))=0\mbox{ for a.e.} \; x \in \O
\]
by the fundamental lemma of the variational method (see~\cite{Willem}). It
follows that $u=0$. So
\[
\varphi(u_n)=\lambda\int_\O\frac{1}{p(x)}|u_n|^{p(x)}dx+\int_\O G(x,u_n)dx\to \lambda\int_\O\frac{1}{p(x)}|u|^{p(x)}dx+\int_\O G(x,u)dx=0.
\]
Hence, we see that for
$t_0=\frac ab$ we have
\[
J(u_n)=a\int_\O\frac{1}{p(x)}|\nabla u_n|^{p(x)}dx-\frac b2\left(\int_\O\frac{1}{p(x)}|\nabla u_n| ^{p(x)}dx\right)^2-\lambda\int_{\O}\frac{1}{p(x)}|u_n|^{p(x)}dx- \int_{\O} G(x,u_n)dx\to \frac{a^2}{2b}.
\]
 This is a contradiction since $J(u_n)\to c<\frac{a^2}{2b}$, then $a-b\int_\O\frac{1}{p(x)}|\nabla u_n|^{p(x)}dx\to 0$ is not true and similarly to Subcase $1$, we have that
\begin{equation}
\{a-b\int_\O\frac{1}{p(x)}|\nabla u_n|^{p(x)}dx\to 0\}  \mbox{ is
bounded}.
\end{equation}
So, it follows from  the two cases above that
\[
\int_\O| \nabla u_n| ^{p(x)-2}\nabla u_n (\nabla u_n-\nabla u)dx\to 0.
\]

Invoking the $S_+$ condition (see Lemma \ref{s+}), we can now deduce that $\| u_n\| \to\| u\| $ as $n\to \infty$, which means
that $J$ satisfies the $(PS)_c$ condition.\qed

\begin{rem}
The $(PS)_c$ condition is not satisfied for $c>\frac{a^2}{2b}$.

 Indeed,
\[
J(u)\leq a\int_\O\frac{1}{p(x)}| \nabla u| ^{p(x)}dx-\frac b2\left(\int_\O\frac{1}{p(x)}| \nabla u| ^{p(x)}dx\right)^2\leq \frac{a^2}{2b}
\]
and so if $\{u_n\}$ is a $(PS)_c$ sequence of $J$, then we have $c\leq \frac{a^2}{2b}$, which is a contradiction.
\end{rem}
\section{Proof of Theorem \ref{theo1.1}}
\label{section4}
To verify
the conditions of the Mountain Pass theorem (see e.g.,~\cite{Willem}),
we first need to prove two lemmas.
\begin{lem}
\label{lemme2} Assume that $g$ satisfies \eqref{hypot:g1} and \eqref{hypot:g2}. Then there exist $\rho > 0$ and
$\alpha > 0$ such that $J(u) \geq \alpha > 0$, for any $u \in W_0^{1,p(x)}(\O)$ with $\| u\| =\rho$.
\end{lem}
{\bf Proof.}
\\
$\bullet$ For $\lambda\leq0$. By assumptions \eqref{hypot:g1} and \eqref{hypot:g2}, we have
\begin{equation}
\label{condg0}
| G(x,u)| \leq\frac{\epsilon}{p(x)}| u| ^{p(x)}+\frac{C_\epsilon}{q(x)}| u| ^{q(x)}.
\end{equation}
Let $\epsilon=\frac18 a\lambda_{p(x)}$  and $u\in W_0^{1,p(x)}(\O)$ be such that $\| u\| =\rho\in(0,1)$. By considering  Lemma \ref{lemmaineq},  Proposition \ref{sobolev} and \eqref{cond q}, we can deduce that
\begin{eqnarray*}
\label{eq0}
J(u)&=&a\int_\O\frac{1}{p(x)}| \nabla u| ^{p(x)}dx-\frac b2\left(\int_\O\frac{1}{p(x)}| \nabla u| ^{p(x)}dx\right)^2-\lambda\int_{\O}\frac{1}{p(x)}|  u| ^{p(x)}dx- \int_{\O} G(x,u)dx\nonumber\\
&\geq&a\int_\O\frac{1}{p(x)}| \nabla u| ^{p(x)}dx-\frac b2\left(\int_\O\frac{1}{p(x)}| \nabla u| ^{p(x)}dx\right)^2-\epsilon \int_{\O} \frac{| u| ^{p(x)}}{p(x)}dx- C_\epsilon \int_{\O} \frac{| u| ^{q(x)}}{q(x)}dx.\nonumber\\
&\geq& (a-\frac{\epsilon}{\lambda_{p(x)}}) \int_{\O}\frac{1}{p(x)}|  \nabla u| ^{p(x)}dx -\frac b2\left(\int_\O\frac{1}{p(x)}| \nabla u| ^{p(x)}dx\right)^2-\frac{CC_\epsilon}{q^-} \int_\O| \nabla u| ^{q(x)}dx
\nonumber \\
&\geq& \frac{1}{p^+}(a-\frac{\epsilon}{\lambda_{p(x)}})\rho_{p(x)}(\nabla u) -\frac{b}{2p^{-^2}}(\rho_{p(x)}(\nabla u))^2-\frac{CC_\epsilon}{q^-}\rho_{q(x)}(\nabla u)\nonumber\\
&\geq& \frac{1}{p^+}(a-\frac{\epsilon}{\lambda_{p(x)}}) \| u\| ^{p^+}-\frac{b}{2p^{{-^2}}}\| u\| ^{2p^-}-\frac{CC_\epsilon}{q^-}\| u\| ^{q^-}\\
&\geq& \left(\frac{7a}{8p^+}-\frac{b}{2p^{-^2}}\| u\| ^{2p^--p^+}-\frac{CC_\epsilon}{q^-}\| u\| ^{q^--p^+}\right)\| u\| ^{p^+}.
\end{eqnarray*}
 We can choose $\rho$ sufficiently small (i.e.~$\rho$ is such that $\frac{7a}{8p^+}-\frac{b}{2p^{-^2}}\rho^{2p^--p^+}-\frac{CC_\epsilon}{q^-}\rho^{q^--p^+}>0$), so that
\[
I(u)\geq \rho^{p^+}(\frac{7a}{8p^+}-\frac{b}{2p^{-^2}}\rho^{2p^--p^+}-\frac{CC_\epsilon}{q^-}\rho^{q^--p^+})\eqcolon \alpha>0.
\]
$\bullet$ For $\lambda>0$. Let $\epsilon>0$ be small enough so that $\frac{1}{2p^+}(a-\frac{\lambda}{\lambda_{p(x)}})=\frac{\epsilon}{\lambda_{p(x)}p^-}$.

 Let $\rho\in(0,1)$ and $u\in W_0^{1,p(x)}(\O)$ be such that $\| u\| =\rho$. By considering  Lemma \ref{lemmaineq}, Proposition \ref{sobolev},  \eqref{cond q}, and \eqref{condg0}, we can deduce that
\begin{eqnarray*}
\label{eq0}
J(u)&=&a\int_\O\frac{1}{p(x)}| \nabla u| ^{p(x)}dx-\frac b2\left(\int_\O\frac{1}{p(x)}| \nabla u| ^{p(x)}dx\right)^2-\lambda\int_{\O}\frac{1}{p(x)}|  u| ^{p(x)}dx- \int_{\O} G(x,u)dx\nonumber\\
&\geq&a\int_\O\frac{1}{p(x)}| \nabla u| ^{p(x)}dx-\frac b2\left(\int_\O\frac{1}{p(x)}| \nabla u| ^{p(x)}dx\right)^2-\frac{\lambda}{\lambda_{p(x)}}\left(\int_{\O}\frac{1}{p(x)}|  \nabla u| ^{p(x)}dx\right)\nonumber\\&&-\epsilon \int_{\O} \frac{| u| ^{p(x)}}{p(x)}dx-C_\epsilon \int_{\O} \frac{| u| ^{q(x)}}{q(x)}dx\nonumber\\
&\geq& (a-\frac{\lambda}{\lambda_{p(x)}})\int_{\O}\frac{1}{p(x)}|  \nabla u| ^{p(x)}dx -\frac b2\left(\int_\O\frac{1}{p(x)}| \nabla u| ^{p(x)}dx\right)^2\nonumber\\&&-\frac{\epsilon}{\lambda_{p(x)}}\int_\O\frac{1}{p(x)}| \nabla u| ^{p(x)}dx -\frac{CC_\epsilon}{q^-} \int_\O| \nabla u| ^{q(x)}dx
\nonumber \\
&\geq&
(\frac{1}{p^+}(a-\frac{\lambda}{\lambda_{p(x)}})-\frac{\epsilon}{\lambda_{p(x)}p^-}) \rho_{p(x)}(\nabla u) -\frac{b}{2p^{-^2}}(\rho_{p(x)}(\nabla u))^2-\frac{CC_\epsilon}{q^-}\rho_{q(x)}(\nabla u)\nonumber\\
&\geq& (\frac{1}{p^+}(a-\frac{\lambda}{\lambda_{p(x)}})-\frac{\epsilon}{\lambda_{p(x)}p^-}) \| u\| ^{p^+}-\frac{b}{2p^{{-^2}}}\| u\| ^{2p^-}-\frac{CC_\epsilon}{q^-}\| u\| ^{q^-}\\
&\geq& \left(\frac{1}{2p^+}(a-\frac{\lambda}{\lambda_{p(x)}})-\frac{b}{2p^{-^2}}\| u\| ^{2p^--p^+}-\frac{CC_\epsilon}{q^-}\| u\| ^{q^--p^+}\right)\| u\| ^{p^+}.
\end{eqnarray*}
Set
\begin{equation}
\lambda^*=\frac{qp^{-^2}\lambda_{p(x)} a
-bp^+q^-\rho^{2p^--p^+}
-2CC_\epsilon p^{-^{2}}
\rho^{q^--p^+}}{q^-p^{-^2}} \mbox{ and } \alpha=\lambda^*\rho^{p^+}.
\end{equation}
We can conclude that for any $\lambda \in (0,\lambda^{*})$, there exists $\alpha > 0$ such that for any $u\in W_0^{1,p(x)}(\O)$ with $\| u\| =\rho$ we have  $J(u) \geq \alpha > 0$.\qed

\begin{lem} \label{lemme3}
 Assume that $g$ satisfies \eqref{hypot:g3}. Then there exists $e\in W_0^{1,p(x)}(\O)$ with $\| e\| >\rho$ (where $\rho$ is given by  Lemma \ref{lemme2}) such that $J(e) < 0$.
\end{lem}
{\bf Proof.} In view of \eqref{hypot:g3} we know that for all $A>0$, there exists $C_A>0$ such that
\begin{equation}\label{d1}
G(x,u)\geq A| u| ^{\theta}- C_A,\; \mbox{for all} \; (x,u)\in
\O\times \mathbb{R}.
\end{equation}
Let $\psi\in C_0^\infty(\O)$, $\psi > 0,$ and $t>1$. By \eqref{d1} we have
\begin{align*}
J(t\psi)&=a\int_\O\frac{1}{p(x)}| t\nabla \psi| ^{p(x)}dx-\frac b2\left(\int_\O\frac{1}{p(x)}| t\nabla \psi| ^{p(x)}dx\right)^2-\lambda\int_{\O}\frac{1}{p(x)}|  t\psi| ^{p(x)}dx- \int_{\O} G(x,t\psi)dx\\&\leq a\int_\O\frac{1}{p(x)}| t\nabla \psi| ^{p(x)}dx-\frac b2\left(\int_\O\frac{1}{p(x)}| t\nabla \psi| ^{p(x)}dx\right)^2-\lambda\int_{\O}\frac{1}{p(x)}| t\psi| ^{p(x)}dx\\&- At^\theta\int_{\O} | t\psi| ^\theta dx + C_A | \O| \\
&\leq \frac{a t^{p^+}}{p^-}\int_\O| \nabla \psi| ^{p(x)}dx-\frac{bt^{2p^-}}{2p^{+^2}}\left(\int_\O| \nabla \psi| ^{p(x)}dx\right)^2-\frac{\lambda}{p^+}t^{p^-}\int_{\O}| \psi| ^{p(x)}dx- At^\theta\int_{\O} | \psi| ^\theta dx + C_A | \O| .
\end{align*}
Since $\theta>2p^->p^+>p^-$, we obtain $J(t\psi)\to -\infty$ $(t\to +\infty)$. Then for $t>1$ large enough, we
can take $e=t\psi$ so that $\| e\| >\rho$ and $J(e) < 0$.\qed
\\
{\bf Proof of  Theorem \ref{theo1.1}}
\\
By Lemmas \ref{lemme1}--\ref{lemme3} and the fact that $J(0) = 0$,
$J$ satisfies the Mountain Pass theorem (see e.g.,~\cite{Willem}).
Therefore, problem \eqref{prob1} has indeed  a nontrivial weak solution.\qed
\section{The proof of Theorem \ref{multiplicity results}}
\label{section5}
The proof mainly rests on an application of the Fountain theorem. Since $X \coloneq  W_0^{1,p(x)}(\O)$ is a separable and reflexive real Banach space, there exist $\{e_j\}\subset X$ and $\{e^*_j\}\subset X^*$ such that
\[
X=\overline{\mbox{span}\{e_j:j=1,2, \ldots \}},\quad X^*=\overline{\mbox{span}\{e^*_j:j=1,2, \ldots \}}
\]
and
\[
\langle e^*_j,e_j\rangle=
\begin{cases}
                          1, & i=j, \\
                          0, & i\neq j.
\end{cases}
\]
For convenience, we write $X_j = \mbox{span}\{e_j\}$, $Y_k=\oplus_{j=1}^k X_j$, $Z_k=\overline{\oplus_{j=k}^\infty X_j}$.
\def\thetheorem{A}
\begin{theorem}[{Fountain Theorem, see~\cite{Willem}}]\label{Lemma2.3}
\label{fountain}
Suppose that an even functional $\Phi \in C^{1}(X,\R)$  satisfies the $(PS)_c$ condition for every $c>0$, and that
 there is $k_{0}>0$ such that for every $k\geq k_{0}$ there exists
$\rho _{k}>r_{k}>0$ so that the following properties hold:
\begin{itemize}
\item[(i)] $a _{k}=\displaystyle{\max_{u\in Y_{k},\| u\|  =\rho _{k}}}\Phi (u)\leq 0$;

\item[(ii)] $b_{k}=\displaystyle{\inf_{u\in Z_{k},\| u\|  =r_{k}}}\Phi (u)\to +\infty $
 as $k\to \infty$.
\end{itemize}
Then $\Phi $ has a sequence of critical points $\{u_{k}\}$ such that
$\Phi(u_{k})\to +\infty$.
\end{theorem}

\begin{lem}
\label{Beta}
Assume that $\alpha\in C_+(\overline{\O})$, $\alpha(x)<p^*(x)$ for any $x\in \overline{\O}$, and denote
\[
\beta_k=\sup_{u\in Z_k, \Vert u\Vert =1}|u| _{\alpha(x)}.
\]
Then $\lim_{k\to\infty} \beta_k=0$.
\end{lem}
{\bf Proof.} Obviously, $0 < \beta_{k+1} \leq \beta_k$, so $\beta_k \to \beta\geq 0$. Let $u_k\in Z_k$ satisfy
\[
\| u_k\| =1,\quad 0\leq \beta_k-| u_k| _{\alpha(x)}<\frac1k.
\]
Then there exists a subsequence of $\{u_k\}$ (which we still denote by $u_k$) such that $u_k \rightharpoonup u$, and
\[
\langle e^*_j,u\rangle=\lim_{k\to\infty}\langle e^*_j,u_k\rangle=0,~~\mbox{for all} \;  e^*_j,
\]
which implies that $u=0$, and so $u_k \rightharpoonup 0$. Since the embedding from $W_0^{1,p(x)} (\O)$ to $L^{\alpha(x)} (\O)$ is compact, it follows that $u_k\to0$ in $L^{\alpha(x)} (\O)$. Hence, we get $\beta_k\to0$ as $k\to\infty$. Proof of  Lemma \ref{Beta} is thus complete.\qed
\\
{\bf Proof of Theorem \ref{multiplicity results}.} By  Lemma \ref{lemme1}, the functional $J$ satisfies the $(PS)_c$ condition where precisely $c<\frac{a^2}{2b}$. Now
we shall verify that $J$ satisfies the conditions of Theorem \ref{fountain} item by item.

$\mathbf{(i)}$\quad By $(g_3)$, there exist $C_1 > 0$, $M > 0$ such that
\begin{equation}
\label{3.19}
G(x,s)\geq C_1| s| ^\theta,\;\;\mbox{for all} \; | s| \geq M,\;\;\;x\in\O.
\end{equation}
Note that by $(g_1)$,
\begin{eqnarray}
\label{3.20}
| G(x,s)| &\leq& \int_0^1| g(x,zs)s| dz\nonumber\\
&\leq&\int_0^1C(1+| zs| ^{q(x)-1})| s| dz\leq C| s| +C| s| ^{q(x)},\;\;\mbox{for all} \;  (x,t)\in\O\times \R.
\end{eqnarray}
Therefore, if $| s|  \leq M$, there exists $C_2 > 0$ such that
\[
| G(x,s)| \leq | s| (C+C| s| ^{q(x)-1})\leq C_2| s| .
\]
Combining this with \eqref{3.19}, we find
\[
G(x,s)\geq C_1| s| ^\theta-C_2| s| ,\;\;\mbox{for all} \;  (x,t)\in\O\times\R.
\]
For $u\in Y_k$, when $\| u\| > 1$,
\begin{eqnarray*}
J(u)&=&a\int_\O\frac{1}{p(x)}| \nabla u| ^{p(x)}dx-\frac b2\left(\int_\O\frac{1}{p(x)}| \nabla u| ^{p(x)}dx\right)^2-\lambda\int_{\O}\frac{1}{p(x)}|  u| ^{p(x)}dx- \int_{\O} G(x,u)dx\nonumber\\
&\leq&a\int_\O\frac{1}{p(x)}| \nabla u| ^{p(x)}dx-\frac b2\left(\int_\O\frac{1}{p(x)}| \nabla u| ^{p(x)}dx\right)^2-\lambda\int_{\O}\frac{1}{p(x)}|  u| ^{p(x)}dx- C_1\int_{\O} | u| ^\theta dx+C_2\int_{\O} | u| dx.
\end{eqnarray*}
Consequently, because when $\| u\| > 1$,
all norms on the finite-dimensional space $Y_k$,  are equivalent, there is $C_W > 0$ such that
\[
\int_{\O}| u| ^{p(x)}dx\geq C_W\| u\| ^{p^-},~~\int_{\O}| u| ^{\theta}dx\geq C_W\| u\| ^{\theta}~\mbox{ and }~\int_{\O}| u| dx\geq C_W\| u\| .
\]
Hence, we get
\[J(u)\leq\frac{a}{p^-}\| u\| ^{p^+}-\frac{b}{2{p^-}^2}\| u\| ^{2^{p-}}-\frac{\lambda C_W}{p^-}\| u\| ^{p^-}-C_1C_W\| u\| ^{\theta}+C_2C_W\| u\| .
\]
Since $\theta>2p^->p^+>p^-$, it follows that for some $\rho_k=\| u\| >0$ large enough we can deduce that
\[
a_{k}=\max_{u\in Y_{k},\| u\|  =\rho _{k}}J (u)\leq 0.
\]
Hence, condition $(i)$ of Theorem \ref{fountain} holds.

$\mathbf{(ii)}$\quad By $(g_1)$ and $(g_2)$, there exist $C_3,C_4>0$ such that
\[
| G(x,u)| \leq\frac{C_3}{p(x)}| u| ^{p(x)}+\frac{C_4}{q(x)}| u| ^{q(x)}.
\]
By computation, we obtain for any  $u\in Z_k$ with $|u| \leq1$,
\begin{eqnarray*}
J(u)&=&a\int_\O\frac{1}{p(x)}| \nabla u| ^{p(x)}dx-\frac b2\left(\int_\O\frac{1}{p(x)}| \nabla u| ^{p(x)}dx\right)^2-\lambda\int_{\O}\frac{1}{p(x)}|  u| ^{p(x)}dx- \int_{\O} G(x,u)dx\\&\geq&a\int_\O\frac{1}{p(x)}| \nabla u| ^{p(x)}dx-\frac b2\left(\int_\O\frac{1}{p(x)}| \nabla u| ^{p(x)}dx\right)^2-\int_{\O}\frac{\lambda}{p(x)}|  u| ^{p(x)}dx- \int_{\O}\frac{C_3}{p(x)} | u| ^{p(x)}dx\\
&&-\int_{\O}\frac{C_4}{q(x)}| u| ^{q(x)}dx\\
&\geq&\frac{a}{p^+}\| u\| ^{p^+}-\frac{b}{2{p^-}^2}\| u\| ^{p^{2-}}-(\lambda+C_3)\frac{\beta_k^{p^-}}{p^-}\| u\| ^{p^-}-\frac{C_4}{q^-}\beta_k^{q^-}\| u\| ^{q^-}.
\end{eqnarray*}
 Let $\varphi\in Z_k$, $\| \varphi\| = 1$ and $0 < t < 1$. Then it follows that
\begin{eqnarray*}
J(t\varphi)&\geq\frac{a}{p^+}t^{p^+}-\frac{b}{2{p^-}^2}t^{p^{2-}}-(\lambda+C_3)\frac{\beta_k^{p^-}}{p^-}t^{p^-}-\frac{C_4}{q^-}\beta_k^{q^-}t^{q^-}\\
&\geq (\frac{a}{p^+}-\frac{b}{2{p^-}^2})t^{p^{2-}}-(\lambda+C_3)\frac{\beta_k^{p^-}}{p^-}t^{p^-}-\frac{C_4}{q^-}\beta_k^{q^-}t^{q^-}.
\end{eqnarray*}
Conditions  $a\geq b$ and $p^+<2{p^-}^2$ imply that $\frac{a}{p^+}-\frac{b}{2{p^-}^2}=\frac{2{p^-}^2a-bp^+}{2{p^-}^2p^+}>0$.

 Hence, we get
\begin{eqnarray*}
J(t\varphi)&\geq (\frac{2{p^-}^2a-bp^+}{2{p^-}^2p^+})t^{p^{2-}}-\frac{C_4}{q^-}\beta_k^{q^-}t^{q^-}-(\lambda+C_3)\frac{\beta_k^{p^-}}{p^-}t^{p^-}\\
&\geq(\frac{2{p^-}^2a-bp^+}{2{p^-}^2p^+}-\frac{C_4}{q^-}\beta_k^{q^-})t^{q^-}-(\lambda+C_3)\frac{\beta_k^{p^-}}{p^-}t^{p^-}.
\end{eqnarray*}
Choosing $\frac{C_4}{q^-}\beta_k^{q^-}<\frac{2{p^-}^2a-bp^+}{4{p^-}^2p^+}$, we can deduce
\[J(t\varphi)\geq\frac{2{p^-}^2a-bp^+}{4{p^-}^2p^+}t^{q^-}-(\lambda+C_3)\frac{\beta_k^{p^-}}{p^-}t^{p^-}.
\]
Obviously, there exists a large enough $k$ such that
\begin{eqnarray*}
J(t\varphi)&\geq&
t^{p^-}\Bigl(\frac{2{p^-}^2a-bp^+}{4{p^-}^2p^+}t^{q^--p^-}-
(\lambda+C_3)\frac{\beta_k^{p^-}}{p^-}\Bigr).
\end{eqnarray*}
Put $\rho_k\coloneq  \left(\frac{4{p^-}^2p^+}{2{p^-}^2a-bp^+}(\lambda+C_3)\frac{\beta_k^{p^-}}{p^-}\right)^{\frac{1}{q^--p^-}}$. Then for sufficiently
large $k$, $\rho_k< 1$. When $t = \rho_k$, $\rho_k\in Z_k$ with $\| \varphi\|  = 1$, we have $J(t\varphi)\geq 0$. Therefore, condition $(ii)$ of  Theorem \ref{fountain} holds. This completes the proof of  Theorem \ref{multiplicity results}.\qed

\subsection*{\bf Acknowledgments}
 Hamdani would like to express his
deepest gratitude to the Military School of Aeronautical Specialities, Sfax (ESA) for providing an excellent
atmosphere for  work. Mtiri  would like to express his gratitude to the Department of Mathematics, Faculty of Sciences
and Arts, King Khalid University, Muhayil Asir, for supporting his work.
Repov\v{s} was supported by the Slovenian Research Agency grants P1-0292, J1-0831,
N1-0064, N1-0083, and N1-0114. The authors wish to acknowledge thes referees for useful comments and valuable suggestions which have helped to improve the presentation.

\end{document}